\documentclass[12pt,twoside,a4paper]{article}
\usepackage[english]{babel}
\usepackage{amssymb}
\usepackage{hyperref}
\usepackage[T1]{fontenc}
\usepackage{amsmath,epsf}
\usepackage{tikz}
\newtheorem{Th}{Theorem}
\newtheorem{Lem}{Lemma}

\newtheorem{Prop}{Proposition}

\newtheorem{Def}{Definition}
\newcommand{\A}{\mathcal{A}}

\newcommand{\F}{\mathcal{F}}
\newcommand{\R}{\mathbb{R}}
\newcommand{\Z}{\mathbb{Z}}

\newcommand{\M}{\mathcal{M}}

\renewcommand{\L}{\mathbb{L}}

\renewcommand{\P}{\mathbb{P}}
\newcommand{\ds}{\displaystyle}
\newcommand{\norm}[1]{\left\lVert #1 \right\rVert}
\newcommand{\E}[1]{\mathbb E\left(#1\right)}
\newcommand{\abs}[1]{\left|#1\right|}
\newcounter{tictac}

\def\1{\,\rlap{\mbox{\small\rm 1}}\kern.15em 1}
\def\ind#1{\1_{#1}}
\def\build#1_#2^#3{\mathrel{\mathop{\kern 0pt#1}\limits_{#2}^{#3}}}
\def\tend#1#2{\build\hbox to 12mm{\rightarrowfill}_{#1\rightarrow #2}^{a.s.}}

\def\converge#1#2#3{\build\hbox to
15mm{\rightarrowfill}_{#1\rightarrow #2}^{\hbox{\scriptsize #3}}}
\begin{document}

\begin{center}
{\Large
    {\sc
Orthomartingale-coboundary decomposition for stationary random fields
    }
}
\bigskip

Mohamed EL MACHKOURI and Davide GIRAUDO\\

\medskip
Laboratoire de Math\'ematiques Rapha\"el Salem\\
UMR CNRS 6085, Universit\'e de Rouen (France)\\
\emph{mohamed.elmachkouri@univ-rouen.fr}\\
\emph{davide.giraudo@etu.univ-rouen.fr}
\end{center}

\bigskip

\begin{abstract}
\baselineskip=18pt We provide a new projective condition for a stationary real random field indexed by the lattice $\Z^d$ to be well approximated by an orthomartingale in the sense of Cairoli (1969). Our
main result can be viewed as a multidimensional version of the martingale-coboundary decomposition
method which the idea goes back to Gordin (1969). It is a powerfull tool for proving limit theorems or large deviations inequalities for stationary random fields when the corresponding result is valid for orthomartingales. \\
\\
\end{abstract}

\section{Introduction and notations}
In probability theory, a powerfull approach for proving limit theorems for stationary sequences of random variables is to find a way to approximate such sequences by martingales. This idea goes back to Gordin \cite{Gordin69} (see Theorem A below). More precisely, let $(X_k)_{k\in\Z}$ be a sequence of real random variables defined on the probability space $(\Omega,\mathcal{F},\mu)$. We assume that $(X_k)_{k\in\Z}$ is stationary in the sense that its finite-dimensional laws are invariant by translations and we denote by $\nu$ the law of $(X_k)_{k\in\Z}$.  Let $f:\R^\Z\to\R$ be defined by $f(\omega)=\omega_0$ and $T:\R^\Z\to\R^\Z$ by $(T\omega)_k=\omega_{k+1}$ for any $\omega$
in $\R^\Z$ and any $k$ in $\Z$. Then the
sequence $(f\circ T^k)_{k\in\Z}$ defined on the probability space $(\R^\Z, \mathcal{B}(\R^\Z), \nu)$ is
stationary with law $\nu$. So, without loss of generality, we can assume that $X_k=f\circ T^k$ for any $k$
in $\Z$. In 1969, Gordin \cite{Gordin69} introduced a powerfull
method
for proving the central limit theorem (CLT) and the weak invariance principle (WIP) for stationary
sequences of dependent random variables satisfying a projective condition (see
\eqref{projective_condition_dim1} below). In the sequel, for any $p\geqslant 1$ and any $\sigma$-algebra $\mathcal{M}\subset\F$, we denote by $\L^p(\Omega,\mathcal{M},\mu)$ the space of $p$-integrable real random variables defined on $(\Omega, \mathcal{M},\mu)$ and we consider the norm $\norm{.}_p$ defined by $\norm{Z}_p^p=\int_{\Omega}\vert Z(\omega)\vert^p d\mu(\omega)$ for any $Z$ in $\L^p(\Omega,\F,\mu)$. We denote also by $\L^p(\Omega,\F,\mu)\ominus\L^p(\Omega,\mathcal{M},\mu)$ the space of all $Z$ in $\L^p(\Omega,\F,\mu)$ such that $\E{Z\mid\mathcal{M}}=0$ a.s.
\newpage
\hspace{-0.8cm}
\textbf{Theorem A (Gordin, 1969)}\label{decomposition_dim_1}
{\em Let $(\Omega, \F, \mu)$ be a probability space and let $T:\Omega\to\Omega$ be a measurable function
such that $\mu=T\mu$. Let also $p\geqslant 1$ and $\M\subset\F$ be a $\sigma$-algebra such that $\M\subset T^{-1}\M$. If $f$ belongs to $\L^p(\Omega,\mathcal{M},\mu)\ominus\L^p(\Omega,\cap_{i\in\Z}T^{-i}\M,\mu)$ such that}
\begin{equation}\label{projective_condition_dim1}
\sum_{k\geqslant 0}\norm{\E{f\mid T^k\mathcal M}}_p<\infty
\end{equation}
{\em then there exist $m$ in $\L^p(\Omega,\M,\mu)\ominus\L^p(\Omega,T\M,\mu)$ and $g$ in $\L^p(\Omega,T\M,\mu)$ such that}
\begin{equation}\label{decomposition_f_dim_1}
f=m+g-g\circ T.
\end{equation}
The term $g-g\circ T$ in \eqref{decomposition_f_dim_1} is called a coboundary and equation
\eqref{decomposition_f_dim_1} is called the martingale-coboundary decomposition of $f$. Moreover, the
stationary sequence $(m\circ T^i)_{i\in\Z}$ is a martingale-difference sequence with respect to the
filtration $(T^{-i}\M)_{i\in\Z}$ (see Definition \ref{Definition_MD_sequence} below) and for any positive integer $n$,
\begin{equation}\label{approximation_partial_sums}
S_n(f)=S_n(m)+g-g\circ T^n
\end{equation}
where $S_n(h)=\sum_{i=0}^{n-1}h\circ T^i$ for any function $h :\Omega\to \R$. Combining
\eqref{approximation_partial_sums} with the Billingsley-Ibragimov CLT for martingales (see
\cite{Billingsley1961} or \cite{Ibragimov1963}), one obtain the CLT for the stationary sequence
$(f\circ T^k)_{k\in\Z}$ when the projective condition \eqref{projective_condition_dim1} holds.
Similarly, combining \eqref{approximation_partial_sums} with the WIP for martingales (see
\cite{Billingsley1968}), we derive the WIP for the stationary sequence $(f\circ T^k)_{k\in\Z}$. Thus,
Gordin's method provides a sufficient condition for proving limit theorems for stationary sequences
when
such a limit theorem holds for martingale-difference sequences. Our aim in this work is to provide an extension of Theorem A for random fields indexed by the lattice $\Z^d$ where $d$ is a positive integer (see Theorem \ref{decomposition_dim_d}).
%
%
%
%
%
%
\section{Main results}
\begin{Def}\label{Definition_MD_sequence}
We say that a sequence $(X_k)_{k\in\Z}$ of real random variables defined on a probability space $(\Omega,\F,\mu)$ is a
martingale-difference (MD) sequence if there exists a filtration
$(\mathcal{G}_k)_{k\in\Z}$ such that $\mathcal{G}_k\subset\mathcal{G}_{k+1}\subset\F$ and $X_k$ belongs to $\L^1(\Omega,\mathcal{G}_k,\mu)\ominus\L^1(\Omega,\mathcal{G}_{k-1},\mu)$ for any $k$ in $\Z$.
\end{Def}
The concept of MD sequences can be extended to the random field setting. One can refer for example to Basu and Dorea \cite{Bas-Dor} or Nahapetian \cite{Nahapetian-Petrosian} where
MD random fields are defined in two differents ways and limit theorems are obtained.
In this paper, we are interested by orthomartingale-difference random fields in the sense of Cairoli
\cite{Cairoli1969}. A good introduction to this concept is done in the book by Khoshnevisan
\cite{Khoshnevisan2002}. Let $d$ be a positive integer. We denote by $\langle d\rangle$ the set
$\{1,...,d\}$. For any $s=(s_1,...,s_d)$ and any $t=(t_1,...,t_d)$ in $\Z^d$, we write $s\preceq t$ (resp. $s\prec t$, $s\succeq t$ and $s\succ t$) if and only if $s_k\leqslant t_k$ (resp. $s_k< t_k$, $s_k\geqslant t_k$ and $s_k>t_k$) for any $k$ in $\langle d\rangle$ and  we denote also $s\wedge t=(s_1\wedge
t_1,...,s_d\wedge t_d)$.
\begin{Def}\label{F4_property}
Let $(\Omega, \F,\mu)$ be a probability space. A family $(\mathcal{G}_i)_{i\in\Z^d}$ of $\sigma$-algebras  is a commuting
filtration if $\mathcal{G}_i\subset\mathcal{G}_j\subset\F$ for any $i$ and $j$ in $\Z^d$ such that
$i\preceq j$ and
$$
\E{\E{Z\mid\mathcal{G}_s}\mid\mathcal{G}_t}
=\E{Z\mid\mathcal{G}_{s\wedge t}}\quad\textrm{a.s.}
$$
for any $s$ and $t$ in $\Z^d$ and any bounded random variable $Z$.
\end{Def}
Definition \ref{F4_property} is known as the ``F4 condition''.
\begin{Def}
Let $(\Omega, \F,\mu)$ be a probability space. A random field $(X_k)_{k\in\Z^d}$ is an orthomartingale-difference (OMD) random field if there exists a commuting filtration $(\mathcal{G}_i)_{i\in\Z^d}$ such that $X_k$ belongs to $\L^1(\Omega,\mathcal{G}_k,\mu)\ominus\L^1(\Omega,\mathcal{G}_{l},\mu)$ for any $l\nsucceq k$ and $k$ in $\Z^d$.
\end{Def}
\textbf{Remark 1.} Let $k$ be fixed in $\Z^d$ and $S_k=\sum_{0\prec i\preceq k}X_i$ where $(X_i)_{i\in\Z^d}$
is an OMD
random field with respect to a commuting filtration $(\mathcal{G}_i)_{i\in\Z^d}$. Then $S_k$ belongs to $\L^1(\Omega,\mathcal{G}_k,\mu)$ and $\E{S_k\mid\mathcal{G}_l}=S_l$ for any $l\preceq k$. We say that $(S_k)_{k\in\Z^d}$ is an orthomartingale (OM) random field.\\
\\
Arguing as above, without loss of generality, every stationary real random field $(X_k)_{k\in\Z^d}$ can
be written as $(f\circ T^k)_{k\in\Z^d}$ where $f:\Omega\to\R$ is a measurable function and for any $k$
in $\Z^d$, $T^k:\Omega\to\Omega$ is a measure-preserving operator satisfying  $T^i\circ T^j=T^{i+j}$
for any $i$ and $j$ in $\Z^d$. For any $s$ in $\langle d\rangle$, we denote $T_s=T^{e_s}$ where
$e_s=(e_s^{(1)},\dots,e_s^{(d)})$ is the unique element of $\Z^d$ such that $e_s^{(s)}=1$ and
$e_s^{(i)}=0$ for any $i$ in $\langle d\rangle\backslash\{s\}$ and $U_s$ is the operator defined  by
$U_sh=h\circ T_s$ for any function $h:\Omega\to\R$. We define also $U_J$ as the product operator
$\Pi_{s\in J}U_s$ for any $\emptyset\subsetneq J\subset \langle d\rangle$ and we write simply $U$ for
$U_{\langle d\rangle}=U_1\circ U_2\circ ...\circ U_d$. For any $\emptyset\subsetneq J\subset \langle
d\rangle$, we denote also by $\vert J\vert$ the number of elements in $J$ and by $J^c$ the set $\langle
d\rangle\backslash J$. The main result of this paper is the following.
\begin{Th}\label{decomposition_dim_d}
Let $(\Omega,\F,\mu)$ be a probability space and let $T^l:\Omega\to\Omega$ be a measure-preserving operator for any $l$ in $\Z^d$ such that $T^i\circ T^j=T^{i+j}$ for any $i$ and $j$ in $\Z^d$. Let $p\geqslant 1$ and let $\mathcal{M}\subset \F$ be a $\sigma$-algebra such that $(T^{-i}\M)_{i\in\Z^d}$
is a commuting filtration. If $f$ belongs to $\L^p(\Omega,\M,\mu)\ominus\L^p(\Omega,\cap_{k\geqslant 1}T_s^k\M,\mu)$ and
\begin{equation}\label{proj_criterion_dim_d}
\sum_{k\geqslant 1} k^{d-1}\norm{\E{f\mid T^k_s\mathcal M}}_p<\infty
\end{equation}
for any $s$ in $\langle d\rangle$ then $f$ admits the decomposition
\begin{equation}\label{decomposition_f_dim_d}
f=m+\sum_{\emptyset\subsetneq J\subsetneq \langle d\rangle}\prod_{s\in J}(I-U_s)m_J+\prod_{s=1}^d(I-U_s)g,
\end{equation}
where $m$, $g$ and $(m_J)_{\emptyset\subsetneq J\subsetneq \langle d\rangle}$ belong to $\L^p(\Omega,\M,\mu)$  and $(U^{i}m)_{i\in \Z^d}$ and
$\left(U_{J^c}^{i}m_J\right)_{i\in\Z^{d-\abs J}}$ for all $\emptyset\subsetneq J\subsetneq \langle
d\rangle$ are OMD random fields.
\end{Th}
\textbf{Remark 2}. If $d=1$ then Theorem \ref{decomposition_dim_d} reduces to Theorem A. If $d=2$ then \eqref{decomposition_f_dim_d} reduces to
$$
f=m+(I-U_1)m_1+(I-U_2)m_2+(I-U_1)(I-U_2)g,
$$
where $m$, $m_1$, $m_2$ and $g$ belong to $\L^p(\Omega,\M,\mu)$ such that $(U^im)_{i\in\Z}$ is an OMD random field and $(U_2^k m_1)_{k\in\Z}$ and $(U_1^k m_2)_{k\in\Z}$ are MD
sequences. If $d=3$ then \eqref{decomposition_f_dim_d} becomes
\begin{align*}
f=m&+(I-U_1)m_1+(I-U_2)m_2+(I-U_3)m_3\\
&+(I-U_1)(I-U_2)m_{\{1,2\}}+(I-U_1)(I-U_3)m_{\{1,3\}}+(I-U_2)(I-U_3)m_{\{2,3\}}\\
&+(I-U_1)(I-U_2)(I-U_3)g
\end{align*}
where $m$, $m_1$, $m_2$, $m_3$, $m_{\{1,2\}}$, $m_{\{1,3\}}$, $m_{\{2,3\}}$ and $g$ belong to $\L^p(\Omega,\M,\mu)$ such that $(U^im)_{i\in\Z^3}$, $(U_{\{2,3\}}^im_1)_{i\in\Z^2}$,
$(U_{\{1,3\}}^im_2)_{i\in\Z^2}$ and $(U_{\{1,2\}}^im_3)_{i\in\Z^2}$ are OMD random fields and $(U_1^k
m_{\{2,3\}})_{k\in\Z}$, $(U_2^k m_{\{1,3\}})_{k\in\Z}$ and $(U_3^k m_{\{1,2\}})_{k\in\Z}$ are MD sequences.\\
\\
\textbf{Remark 3}. A decomposition similar to \eqref{decomposition_f_dim_d}  was obtained by Gordin \cite{Gordin2009} but with reversed orthomartingales and  under an assumption on the so-called Poisson equation.
\begin{Prop}\label{Moment_inequalities}
Let $(X_i)_{i\in\Z^d}$ be an OMD random field. There exists a positive constant $\kappa$ such that for any $p>1$ and any $n$ in $\Z_+^d$,
\begin{equation}\label{Moment_inequality_OMD}
\norm{\sum_{0\preceq k\preceq n}X_k}_{p}\leqslant\kappa p^{d/2}\left(\sum_{0\preceq k\preceq n}\norm{X_k}_p^2\right)^{1/2}
\end{equation}
and the constant $p^{d/2}$ in \eqref{Moment_inequality_OMD} is optimal in the following sense: there exists a stationary OMD random field $(Z_k)_{k\in\Z^d}$ with $\norm{Z_0}_{\infty}=1$ and a positive constant $\kappa$ such that for any $p>1$
\begin{equation}\label{borne_inf_constants}
\inf\left\{C>0\,\,;\,\,\norm{\sum_{0\preceq k\preceq n}Z_k}_{p}\leqslant C\left(\sum_{0\preceq k\preceq n}\norm{Z_k}_p^2\right)^{1/2}\forall{n\in\Z_+^d}\right\}
\geqslant \kappa p^{d/2}.
\end{equation}
\end{Prop}
\textbf{Remark 4}. A Young function $\psi$ is a real convex nondecreasingfunction defined on $\R^{+}$ which satisfies $\lim_{t\to\infty}\psi(t)=\infty$ and $\psi(0)=0$. We denote by $\L_{\psi}(\Omega,\F,\mu)$ the Orlicz space associated to the Young function $\psi$, that is the space of real random
variables $Z$ defined on $(\Omega, \F, \mu)$
such that $\E{\psi(\vert Z\vert/c)}<\infty$ for some $c>0$. The
Orlicz space $\L_{\psi}(\Omega,\F,\mu)$ equipped with the so-called Luxemburg
norm $\| . \|_{\psi}$ defined for any real random
variable $Z$ by $\| Z\|_{\psi}=\inf\{\,c>0\,;\,E[\psi(\vert
Z\vert/c)]\leqslant 1\,\}$ is a Banach space. For any $p\geqslant 1$, if $\varphi_p$ is the function defined by $\varphi_p(x)=x^p$ for any nonegative real $x$ then $\varphi_p$ is a Young function and the Orlicz space $\L_{\varphi_p}(\Omega, \F,\mu)$ reduces to $\L^p(\Omega,\F,\mu)$ equipped with the norm $\norm{.}_p$ defined by $\norm{Z}_p^p=\int_{\Omega}\vert Z(\omega)\vert^p d\mu(\omega)$ for any $Z$ in $\L^p(\Omega,\F,\mu)$. For more about Young functions and Orlicz
spaces one can refer to Krasnosel'skii and Rutickii \cite{K-R}. Combining Lemma $4$ in \cite{Elmachkouri--Volny--Wu2013} and Inequality \eqref{Moment_inequality_OMD}, we obtain Kahane-Khintchine inequalities: for any $0<q<2/d$, there exists a positive constant $\kappa_q$ such that for any $n\succeq 0$ in $\Z^d$,
\begin{equation}\label{KK_inequalities_OMD}
\norm{\sum_{0\preceq k\preceq n}X_k}_{\psi_q}\leqslant\kappa_q\left(\sum_{0\preceq k\preceq n}\norm{X_k}_{\psi_{\beta(q)}}^2\right)^{1/2}
\end{equation}
where $\beta(q)=2q/(2-dq)$ and $\psi_\alpha$ is the Young function defined for any $x\in\R^{+}$ by
$$
\psi_{\alpha}(x)=\exp((x+h_{\alpha})^{\alpha})-\exp(h_{\alpha}^{\alpha})
\quad\textrm{with}\quad h_{\alpha}=((1-\alpha)/\alpha)^{1/\alpha}
\ind{\{0<\alpha<1\}}
$$
for any real $\alpha>0$. Moreover, \eqref{KK_inequalities_OMD} still hold for $q=2/d$ if the random variables $(X_k)_{k\in\Z^d}$ are assumed to be uniformly bounded. Finally, using Markov inequality and the definition of the Luxembourg norm $\norm{.}_{\psi_{q}}$, we derive the following large deviations inequalities: for any $0<q<2/d$, there exists a positive constant $\kappa_q$ such that for any $n\succeq 0$ in $\Z^d$ and any $x>0$,
$$
\mu\left(\left\vert\sum_{0\preceq k\preceq n}X_k\right\vert>x\right)\leqslant(1+e^{h_q^q})
\exp\left(-\left(\frac{x}{\kappa_q\sqrt{\sum_{0\preceq k\preceq n}\norm{X_k}_{\psi_{\beta(q)}}^2}}+h_q\right)^q\right).
$$
Again, the above exponential inequality still hold for $q=2/d$ when the random variables $(X_k)_{k\in\Z^d}$ are uniformly bounded.\\
\\
Combining Proposition \ref{Moment_inequalities} and Theorem \ref{decomposition_dim_d}, we obtain the following result.
\begin{Prop}\label{Moment_inequality_under_proj_condition}
Let $(X_i)_{i\in\Z^d}$ be a stationary real random field and $(\F_i)_{i\in\Z^d}$ be a commuting filtration such that each $X_i$ is $\F_i$-measurable. If there exists $p>1$ such that $X_0$ belongs to $\L^p(\Omega,\F_0,\mu)\ominus\L^p(\Omega,\cap_{k\geqslant 1}\F_{k,s},\mu)$ and
\begin{equation}\label{criterion_IP_dim_d_norm_p}
\sum_{k\geqslant 1} k^{d-1}\norm{\E{X_0\mid \F_{k,s}}}_p<\infty
\end{equation}
for any $s$ in $\langle d\rangle$ where $\F_{k,s}=\vee_{\substack{i=(i_1,..,i_d)\in\Z^d\\ i_s\leqslant -k}}\F_{i}$ then \eqref{Moment_inequality_OMD} still holds.
\end{Prop}
Now, we are able to investigate the WIP for random fields. Let $(X_i)_{i\in\Z^d}$ be a stationary real random field defined on a probability space $(\Omega,\F,\mu)$. Let also $\A$ be a collection of Borel subsets of $[0,1]^{d}$ and consider the process $\{S_{n}(A)\,;\,A\in\A\}$ defined by
\begin{equation}\label{process}
S_{n}(A)=\ds{\sum_{i\in\langle n\rangle^d}}\,\lambda(nA\cap R_{i})X_{i}\end{equation}
where $R_{i}=]i_{1}-1,i_{1}]\times...\times]i_{d}-1,i_{d}]$ is the unit cube with upper corner at $i=(i_1,..,i_d)$ in $\langle n\rangle^d$ and $\lambda$ is the Lebesgue measure on $\R^{d}$. The collection $\A$ is equipped with the
pseudo-metric $\rho$ defined by $\rho(A,B)=\sqrt{\lambda(A\Delta B)}$ for any $A$ and $B$ in $\A$. Let $\varepsilon>0$ and let $H(\A,\rho,\varepsilon)$ be the logarithm of the smallest number
$N(\A,\rho,\varepsilon)$ of open balls of radius $\varepsilon$ with respect to $\rho$ which form a covering of $\A$. The function
$H(\A, \rho, .)$ is called the metric entropy of the class $\A$ and allows us to control the size of the collection $\A$. Let $(\mathcal{C}(\A),\|.\|_{\A})$ be the Banach space of continuous real functions on $\A$ equipped with the uniform norm
$\|.\|_{\A}$ defined by $\| f\|_{\A}=\sup_{A\in\A}\vert f(A)\vert$. A standard Brownian motion indexed by $\A$ is a mean zero Gaussian
process $W$ with sample paths in $\mathcal{C}(\A)$ and Cov$(W(A),W(B))=\lambda(A\cap B)$. From Dudley \cite{Dudley} we know that such a process is well defined if $
\int_{0}^{1}\sqrt{H(\A,\rho,\varepsilon)}\,d\varepsilon<\infty$.
Following \cite{van-der-Vaart-Wellner}, we recall the definition of Vapnik-Chervonenkis classes ($VC$-classes) of sets: let $\mathcal{C}$ be a collection
of subsets of a set $\mathcal{X}$. An arbitrary set of $n$ points $F_n:=\{x_1,...,x_n\}$ possesses $2^n$ subsets.
Say that $\mathcal{C}$ {\em picks out} a certain subset from $F_n$ if this can be formed as a set of the form $C\cap F_n$
for a $C$ in $\mathcal{C}$. The collection $\mathcal{C}$ is said to {\em shatter} $F_n$ if each of its $2^n$ subsets can be picked out in this manner. The {\em VC-index}
$V(\mathcal{C})$ of the class $\mathcal{C}$ is the smallest $n$ for which no set of size $n$ is shattered by $\mathcal{C}$. Clearly,
the more refined $\mathcal{C}$ is, the larger is its index. Formally, we have
$$
V(\mathcal{C})=\inf\left\{n\,;\,\max_{x_1,...,x_n}\Delta_n(\mathcal{C},x_1,...,x_n)<2^n\right\}
$$
where $\Delta_n(\mathcal{C},x_1,...,x_n)=\#\left\{C\cap\{x_1,...,x_n\}\,;\,C\in\mathcal{C}\right\}$. Two classical examples of
$VC$-classes are the collection $\mathcal{Q}_d=\left\{[0,t]\,;\,t\in[0,1]^{d}\right\}$ and $\mathcal{Q}^{'}_d=\left\{[s,t]\,;\,s,t\in[0,1]^{d},\,s\preceq t\right\}$
with index $d+1$ and $2d+1$ respectively.\\
\\
In the sequel, since the CLT does not hold for general OMD random fields (see \cite{Wang--Woodroofe2011b}, example 1, page 12), we restrict ourselves to the case of a commuting filtration generated by iid random variables.
\begin{Th}\label{IP_OMD}
Let $(\varepsilon_j)_{j\in\Z^d}$ be an iid real random field defined on a probability space $(\Omega,\F,\mu)$. Denote by $(\F_i)_{i\in\Z^d}$ the commuting filtration where $\F_i$ is the $\sigma$-algebra generated by $\varepsilon_j$ for $j\preceq i$ and $i$ in $\Z^d$ and consider the $\sigma$-algebra $\F_{k,s}$ defined in Proposition \ref{Moment_inequality_under_proj_condition} for any positive integer $k$ and any $s$ in $\langle d\rangle$. Let $(X_i)_{i\in\Z^d}$ be a stationary real random field and $\A$ be a VC-class of regular Borel subsets of $[0,1]^d$ with index $V$. Assume that there exists $p>2(V-1)$ such that $X_0$ belongs to $\mathbb{L}^p(\Omega,\F_0,\mu)\ominus\mathbb{L}^p(\Omega,\cap_{k\geqslant 1}\F_{k,s},\mu)$ for any $s$ in $\langle d\rangle$ and \eqref{criterion_IP_dim_d_norm_p} holds. Then the sequence of processes $\{n^{-d/2}S_n(A)\,;\,A\in\A\}$ converges in distribution in $\mathcal{C}(\A)$
to $\sqrt{\E{X_0^2}} W$ where $W$ is a standard Brownian motion indexed by $\A$.
\end{Th}
If we consider the particular case $\A=\mathcal{Q}_d$ where $\mathcal{Q}_d$ is the class of all quadrants $[0,t]$ for $t$ in $[0,1]^d$ then Theorem \ref{IP_OMD} ensures that the WIP holds for $p$-integrable OMD random fields with $p>2d$. In fact, our next result shows that the WIP still holds for $p=2$ and can be viewed as an extension of the Donsker's invariance principle for iid random variables (see \cite{Donsker1951}). 
\begin{Th}\label{IP_OMD_bis}
Theorem \ref{IP_OMD} still holds with $\A=\mathcal{Q}_d$ and $p=2$.
\end{Th}
\textbf{Remark 5}. El Machkouri et al. \cite{Elmachkouri--Volny--Wu2013} and Wang and Woodroofe
\cite{Wang--Woodroofe2011b} obtained also a WIP for random fields $(X_k)_{k\in\Z^d}$ which can be expressed as a functional of iid real random variables but under the more restrictive condition that $X_0$ belongs to $\L^p(\Omega,\F,\mu)$ with $p>2$. In a recent work, Wang and Voln\'y
\cite{Volny--Wang2014} obtained the WIP for $p=2$ under a multidimensional version of the so-called
Hannan's condition for time series. Their condition is less restrictive than (\ref{criterion_IP_dim_d_norm_p})
but condition (\ref{criterion_IP_dim_d_norm_p}) gives also an orthomartingale approximation for the considered random field which is of independent interest.
\begin{Prop}\label{IP_linear_random_fields}
Let $(\varepsilon_i)_{i\in\Z^d}$ be an iid real random field defined on a probability space $(\Omega,\F,\mu)$ such that $\varepsilon_0$ belongs to $\L^p(\Omega,\F,\mu)$ for some $p\geqslant 2$. Consider the linear random field $(X_k)_{k\in\Z^d}$ defined for any $k$ in $\Z^d$ by $
X_k=\sum_{j\succeq 0}a_j\varepsilon_{k-j}$ where $(a_j)_{j\in\Z^d}$ is a family of real numbers satisfying $\sum_{j\succeq 0}a_j^2<\infty$. Then the condition
\begin{equation}\label{condition_IP_linear_random_field}
\sum_{k\geqslant 1}k^{d-1}\sqrt{\sum_{i\in\Lambda_{k,s}}a_i^2}
\end{equation}
where $\Lambda_{k,s}=\{i=(i_1,..,i_d)\in\Z^d\,;\,i_s\geqslant k\}$ is more restrictive than \eqref{criterion_IP_dim_d_norm_p}.
\end{Prop}
In particular, Proposition \ref{IP_linear_random_fields} ensures that the conclusions of Theorem \ref{IP_OMD} and Theorem \ref{IP_OMD_bis} still hold for linear random fields with iid innovations under assumption \eqref{condition_IP_linear_random_field}. In the other part,
Wang and Woodroofe (\cite{Wang--Woodroofe2011b}, Corollary 1) obtained a WIP for stationary linear
random fields with iid innovations $(\varepsilon_i)_{i\in\Z^d}$ under a weaker condition than (\ref{condition_IP_linear_random_field}) but again with the
additional assumption that $\varepsilon_0$ belongs to $\L^p(\Omega,\F,\mu)$ with $p>2$.\\
\\
%
%
%
%
%
%
We now provide an application of Theorem \ref{decomposition_dim_d} to the WIP in
H\"older spaces. We consider for $0<\beta\leqslant 1$ the space $\mathbb{H}_\beta([0,1]^d)$ as
the space of all continuous functions $g$ for which there exists a constant $K$ such that
for each $s,t\in [0,1]^d$,
\[
|g(s)-g(t)|\leqslant K\lVert s-t\rVert^\beta,
\]
where $\lVert\cdot\rVert$ denotes the Euclidian norm on $\mathbb R^d$. We endow this function
space with the norm $\lVert g\rVert:=|g(0)|+\sup_{s,t\in [0,1]^d,s\neq t}|g(t)-g(s)|
/\lVert t-s\rVert^\beta$. If $(X_k)_{k\in\Z^d}$ is a stationary real random field, we define the partial sum process $\{n^{-d/2}S_n(t);t\in[0,1]^d\}_{n\geqslant 1}$ by 
\begin{equation}\label{Definition_partial_sum_process}
S_n(t)=\sum_{i\in\langle n\rangle^d}\lambda([0,nt]\cap R_i)X_i
\end{equation}
for any positive integer $n$ and any $t$ in $[0,1]^d$ and we recall that $\lambda$ is the Lebesgue measure on $\R^d$ and $R_i=]i_1-1,i_1]\times...\times]i_d-1,i_d]$ is
the unit cube with upper corner $i=(i_1,\dots,i_d)$ in $\langle n\rangle^d$. We consider $(S_n(t))_{t\in [0,1]^d}$  as an
element of $\mathbb{H}_\beta([0,1]^d)$. Our next result provides a sufficient condition
for the weak convergence of $(n^{-d/2}S_n(\cdot))_{n\geqslant 1}$ in this function space.

\begin{Th}\label{IP_Holder}
If the assumptions of Theorem \ref{IP_OMD} hold with $p>4\times (\log_2(4d/(4d-3)))^{-1}$ and $\A=\mathcal{Q}_d$ then the WIP holds in $\mathbb{H}_\gamma([0,1]^d)$ for
each $\gamma<1/2-d/p$.
\end{Th}
\textbf{Remark 6}. In \cite{Suquet-Zemlys-Rackauskas2007}, a necessary and sufficient condition was obtained for iid random fields to satisfy the WIP in H\"older spaces. Our result provides a sufficient condition for stationary real random fields which can be expressed as a functional of iid real random variables.
\section{Proofs}
In this section, the letter $\kappa$ will denote a universal positive constant which the value may change from line to line. The following two lemmas will be useful in the sequel.
\begin{Lem}\label{indep}
Let $(\Omega, \F,\mu))$ be a probability space. For mutually independent sub-$\sigma$-algebras $\mathcal A_1,\mathcal A_2$ and $\mathcal A_3$ of $\F$ and $X$
a random variable in $\L^1(\Omega,\F,\mu)$, we have
\begin{equation}
\E{\E{X\mid \mathcal A_1\vee \mathcal A_2}\mid \mathcal A_2\vee \mathcal A_3}=\E{X\mid \mathcal A_2},
\end{equation}
where $\mathcal A_1\vee\mathcal A_2$ is the $\sigma$-algebra generated by $\mathcal A_1$ and $\mathcal
A_2$.
\end{Lem}
For a proof of Lemma \ref{indep}, one can refer to Proposition 8.1 in \cite{Wang--Woodroofe2011b}.
%
%
\begin{Lem}\label{lemme_cle}
Let $(\Omega,\F,\mu)$ be a probability space and let $T^l:\Omega\to\Omega$ be a measure-preserving operator for any $l$ in $\Z^d$ such that $T^i\circ T^j=T^{i+j}$ for any $i$ and $j$ in $\Z^d$ ($d\geqslant 2$). Let $p\geqslant 1$ be fixed  and let $\mathcal{M}\subset \F$ be a $\sigma$-algebra such that $(T^{-i}\M)_{i\in\Z^d}$
is a commuting filtration. Assume that $F$ belongs to
$\mathbb{L}^p(\Omega, \mathcal{M},\mu)$ and for any $s$ in $\langle d\rangle$,
\begin{equation}\label{assumption_serie_p}
\sum_{k\geqslant 1}k^{d-1}\norm{\mathbb{E}\left(F\mid T_s^k\mathcal M\right)}_p<\infty.
\end{equation}
Then for any $s$ in $\langle d\rangle$, there exist $M_s$ in $\L^p(\Omega, \mathcal{M},\mu)\ominus\L^p(\Omega, T_s\mathcal{M},\mu)$ and
$G_s$ in $\mathbb{L}^p(\Omega, T_s\mathcal{M},\mu)$ such that $F=M_s+(I-U_s)G_s$. Moreover, for any $s$ and $l$ in $\langle d\rangle$, we have
\begin{equation}
\sum_{k\geqslant 1} k^{d-2}\norm{\mathbb{E}\left(M_s\mid T_l^k\mathcal
M\right)}_p<\infty\quad\textrm{and}\quad
\sum_{k\geqslant 1} k^{d-2}\norm{\E{G_s\mid T_l^k\mathcal M}}_p<\infty.
\end{equation}
\end{Lem}
%
%
%
%
{\em Proof of Lemma \ref{lemme_cle}}. The first part of the proposition is well known (see
\cite{Volny93}, Theorem 2). In fact, \eqref{assumption_serie_p} is a sufficient condition for $F$ to
be equal to $M_s+(I-U_s)G_s$ with $M_s$ in $\L^p(\Omega, \mathcal{M},\mu)\ominus\L^p(\Omega, T_s\mathcal{M},\mu)$ and $G_s\in \mathbb{L}^p(\Omega, T_s\mathcal{M},\mu)$
for any $s$ in $\langle d\rangle$. Moreover, $M_s$ and $G_s$ are given by$$
M_s=\sum_{k\geqslant 0} \E{U_s^kF\mid \mathcal M}-\E{U_s^kF\mid T_s\mathcal M}\quad\textrm{and}\quad
G_s=\sum_{k\geqslant 0}\E{U_s^kF\mid T_s\mathcal M}.
$$
Let $r$ be a positive integer and let $s$ and $l$ be fixed in $\langle d\rangle$. We have
$$
\norm{\E{M_s\mid T_l^r\mathcal M}}_p
\leqslant 2\sum_{k\geqslant 0} \norm{\E{F\mid T_{s}^kT_l^r\mathcal M}}_p\leqslant 2r\norm{\E{F\mid
T_l^r\mathcal M}}_p+2\sum_{k\geqslant r}
\norm{\E{F\mid T_s^k\mathcal M}}_p
$$
and consequently
\begin{align*}
\sum_{r\geqslant 1}r^{d-2}\norm{\E{M_s\mid T_l^r\mathcal M}}_p&\leqslant 2\sum_{r\geqslant 1} r^{d-1}\norm{\E{F\mid
T_l^r\mathcal M}}_p+2\sum_{k\geqslant 2}\sum_{r=1}^{k}r^{d-2}\norm{\E{F\mid T_s^k\mathcal M}}_p\\
&\leqslant 2\sum_{r\geqslant 1} r^{d-1}\norm{\E{F\mid T_l^r\mathcal M}}_p+2\sum_{k\geqslant 2} k^{d-1}\norm{\E{F\mid
T_s^r\mathcal M}}_p<\infty.
\end{align*}
Similarly, we have also $\sum_{r\geqslant 1} r^{d-2}\norm{\E{G_s\mid T_l^r\mathcal M}}_p<\infty$. The proof of Lemma \ref{lemme_cle} is complete.\\
\\
%
%
%
{\em Proof of Theorem \ref{decomposition_dim_d}}.
For simplicity, we consider only the case $d=2$ and the case $d=3$. We start with $d=2$. Since $f$ is $\mathcal M$-measurable and
$\sum_{k=1}^{\infty}k\norm{\E{f\mid T_2^k\mathcal M}}_p<\infty$, there exist two functions $m_2$ and
$g_2$ (see \cite{Volny93}, Theorem 2) such that
\begin{equation}\label{decomposition_f}
f=m_2+(I-U_2)g_2,
\end{equation}
where $m_2\in \mathbb{L}^p(\Omega,\M,\mu)\ominus\mathbb{L}^p(\Omega,T_2\M,\mu)$ and $g_2\in \mathbb{L}^p(\Omega,T_2\M,\mu)$. We lay emphasis on that a careful reading of the proof of Theorem 2 in
\cite{Volny93} ensure that $g_2$ is $T_2\M$-measurable when $f$ is $\M$-measurable. So, by Lemma
\ref{lemme_cle}, we have $
\sum_{r\geqslant 1}\norm{\E{m_2\mid T_1^r\mathcal M}}_p
<\infty$. Applying again Theorem 2 in \cite{Volny93}, we obtain
\begin{equation}\label{decomposition_m_seconde}
m_2=m_1+(I-U_1)g_1,
\end{equation}
where $m_1\in \mathbb{L}^p(\Omega,\M,\mu)\ominus\mathbb{L}^p(\Omega,T_1\M,\mu)$ and $g_1\in \mathbb{L}^p(\Omega,T_1\M,\mu)$. Consequently,
\begin{equation}\label{decomposition_f_m1_g1_g2}
f=m_1+(I-U_1)g_1+(I-U_2)g_2.
\end{equation}
Put $m:=m_1-\E{m_1\mid T_2\mathcal M}$ so that $\E{m\mid T_2\mathcal M}=0$. Keeping in mind that
$\E{m_1\mid T_1\mathcal M}=0$ and that $(T^{-i}\mathcal{M})_{i\in\Z^d}$ is a commuting filtration, we derive $
\E{\E{m_1\mid T_2\mathcal M}\mid T_1\mathcal M}=\E{m_1\mid T_1T_2\mathcal M}=0$ and consequently
$\E{m\mid T_1\mathcal M}=0$. That is, $(U^im)_{i\in\Z^2}$ is an OMD random field. Using
(\ref{decomposition_m_seconde}) and $\E{m_2\mid T_2\mathcal M}=0$, we derive
$$
\E{m_1\mid T_2\mathcal M}
=-\E{g_1\mid T_2\mathcal M}+U_1\E{g_1\mid T_1T_2\mathcal M}.
$$
Since $g_1$ is $T_1\mathcal{M}$-measurable and $(T^{-i}\mathcal{M})_{i\in\Z^d}$ is a commuting filtration, we have $\E{g_1\mid
T_1T_2\mathcal M}=\E{g_1\mid T_2\mathcal M}$ and consequently we obtain\begin{equation}\label{terme_correcteur}
\E{m_1\mid T_2\mathcal M}
=-(I-U_1)\E{g_1\mid T_2\mathcal M}.
\end{equation}
Combining (\ref{decomposition_f_m1_g1_g2}) and (\ref{terme_correcteur}), we get
\begin{equation}\label{decomposition_f_temp}
f=m+(I-U_1)\left[g_1-\E{g_1\mid T_2\mathcal{M}}\right]+(I-U_2)g_2
\end{equation}
Now, it suffices to find a decomposition for $g_2$ with respect to $T_1$. In fact, by Lemma
\ref{lemme_cle}, we have $\sum_{r\geqslant 1}\norm{\E{g_2\mid T_1^r\mathcal M}}_p<\infty$. So, we can write
(see \cite{Volny93}, Theorem 2),
\begin{equation}\label{decomposition_g_seconde}
g_2=\overline{m}_1+(I-U_1)\overline{g}_1,
\end{equation}
where $\overline{m}_1\in \mathbb{L}^p(\Omega,\M,\mu)\ominus\mathbb{L}^p(\Omega,T_1\M,\mu)$ and $\overline{g}_1\in \mathbb{L}^p(\Omega,T_1\M,\mu)$. That is, $(U_1^k\overline{m}_1)_{k\in\Z}$ is a MD
sequence. Combining (\ref{decomposition_f_temp}) and (\ref{decomposition_g_seconde}), we derive
$$
f=m+(I-U_1)\left[g_1-\E{g_1\mid T_2\mathcal{M}}\right]+(I-U_2)\overline{m}_1+(I-U_2)(I-
U_1)\overline{g}_1.
$$
The proof of Theorem \ref{decomposition_dim_d} is complete for $d=2$.\\\\
%
%
%
%
In order to convince the reader, we consider now the case $d=3$. Applying again Theorem 2 in
\cite{Volny93}, we decompose $f$ in the following way:
\begin{equation}\label{decomposition_f_en_m3_g3}
f=m_3+(I-U_3)g_3
\end{equation}
where $m_3\in \mathbb{L}^p(\Omega,\M,\mu)\ominus\mathbb{L}^p(\Omega,T_3\M,\mu)$ and $g_3\in \mathbb{L}^p(\Omega,T_3\M,\mu)$. By Lemma \ref{lemme_cle}, we have $\sum_{k\geqslant 1}k\norm{\E{m_3\mid T_2^k\mathcal
M}}_p<\infty$. Since $m_3$ is $\mathcal M$-measurable, we obtain
\begin{equation}\label{decomposition_m3_en_m2_g2}
m_3=m_2+(I-U_2)g_2
\end{equation}
where $m_2\in \mathbb{L}^p(\Omega,\M,\mu)\ominus\mathbb{L}^p(\Omega,T_2\M,\mu)$ and $g_2\in \mathbb{L}^p(\Omega,T_2\M,\mu)$. By Lemma \ref{lemme_cle}, we have also $\sum_{k\geqslant 1}\norm{\E{m_2\mid T_1^k\mathcal
M}}_p<\infty$. Since $m_2$ is $\mathcal M$-measurable, we have also
\begin{equation}\label{decomposition_m2_en_m1_g1}
m_2=m_1+(I-U_1)g_1
\end{equation}
where $m_1\in \mathbb{L}^p(\Omega,\M,\mu)\ominus\mathbb{L}^p(\Omega,T_1\M,\mu)$ and $g_1\in \mathbb{L}^p(\Omega,T_1\M,\mu)$. Combining (\ref{decomposition_f_en_m3_g3}) ,(\ref{decomposition_m3_en_m2_g2}) and
(\ref{decomposition_m2_en_m1_g1}), we obtain
\begin{equation}\label{decomposition_f_en_m1_g1_g2_g3}
f=m_1+(I-U_1)g_1+(I-U_2)g_2+(I-U_3)g_3.
\end{equation}
Now, we define
\begin{equation}\label{definition_m}
m:=m_1-\E{m_1\mid T_2\mathcal M}-\E{m_1\mid T_3\mathcal M}+\E{m_1\mid T_2T_3\mathcal M}.
\end{equation}
Since $\E{m_1\mid T_1\mathcal M}=0$ and $(T^{-i}\mathcal{M})_{i\in\Z^3}$ is a commuting filtration, $m$ satisfies
$$
\E{m\mid T_1\mathcal M}=\E{m\mid T_2\mathcal M}=
\E{m\mid T_3\mathcal M}=0.
$$
So, $(U^im)_{i\in\Z^3}$ is an OMD random field. Moreover, using (\ref{decomposition_m2_en_m1_g1}),
$$
\E{m_1\mid T_2\mathcal M}=\E{m_2-(I-U_1)g_1\mid T_2\mathcal M}=-(I-U_1)\E{g_1\mid T_1T_2\mathcal M}.
$$
Since again $g_1$ is $T_1\mathcal{M}$-measurable and $(T^{-i}\mathcal{M})_{i\in\Z^3}$ is a commuting filtration, we have $\E{g_1\mid T_1T_2\mathcal
M}=\E{g_1\mid T_2\mathcal M}$ and consequently,
\begin{equation}\label{esperance_m1_sachant_T2M}
\E{m_1\mid T_2\mathcal M}=-(I-U_1)\E{g_1\mid T_2\mathcal M}.
\end{equation}
Similarly,
\begin{equation}\label{esperance_m1_sachant_T3M}
\E{m_1\mid T_3\mathcal{M}}=-(I-U_1)\E{g_1\mid T_3\mathcal{M}}-(I-U_2)\E{g_2\mid T_3\mathcal{M}}
\end{equation}
and
\begin{equation}\label{esperance_m1_sachant_T2T3M}
\E{m_1\mid T_2T_3\mathcal{M}}=-(I-U_1)\E{g_1\mid T_2T_3\mathcal{M}}.
\end{equation}
Combining (\ref{decomposition_f_en_m1_g1_g2_g3}), (\ref{definition_m}),
(\ref{esperance_m1_sachant_T2M}), (\ref{esperance_m1_sachant_T3M}) and
(\ref{esperance_m1_sachant_T2T3M}), we obtain
\begin{align*}
f=m&+(I-U_1)[g_1-\E{g_1\mid T_2\M}-\E{g_1\mid T_3\M} +\E{g_1\mid T_2T_3\M}]\\
&+(I-U_2)[g_2-\E{g_2\mid T_3\M}]+(I-U_3)g_3.
\end{align*}
By (\ref{decomposition_f_en_m3_g3}) and Lemma \ref{lemme_cle}, we have $\sum_{k\geqslant 1}k\norm{\E{g_3\mid
T_1^k\mathcal M}}_p<\infty$. Since $g_3$ is $T_3\M$-measurable, we derive from Theorem 2 in
\cite{Volny93} that
\begin{equation}\label{decomposition_g3_en_m1barre_g1barre}
g_3=\overline{m}_1+(I-U_1)\overline{g}_1
\end{equation}
where $\overline{m}_1\in \mathbb{L}^p(\Omega, T_3\mathcal M,\mu)\ominus\mathbb{L}^p(\Omega, T_1T_3\mathcal M,\mu)$ and $\overline{g}_1\in
\mathbb{L}^p(\Omega, T_1T_3\mathcal M,\mu)$. By Lemma
\ref{lemme_cle}, we have $\sum_{k\geqslant 1}\norm{\E{\overline{m}_1\mid T_2^k\mathcal M}}_p<\infty$ and
consequently
\begin{equation}\label{decomposition_m1barre_en_m2barre_g2barre}
\overline{m}_1=\overline{m}_2+(I-U_2)\overline{g}_2
\end{equation}
where $\overline{m}_2\in \mathbb{L}^p(\Omega, T_3\mathcal M,\mu)\ominus\mathbb{L}^p(\Omega, T_2T_3\mathcal M,\mu)$ and $\overline{g}_2\in
\mathbb{L}^p(\Omega, T_2T_3\mathcal M,\mu)$.
Denoting
\begin{equation}\label{definition_mbarre}
\overline{m}:=\overline{m}_2-\E{\overline{m}_2\mid T_1T_3\M},
\end{equation}
and keeping in mind that $(T^{-i}\mathcal{M})_{i\in\Z^d}$ is a commuting filtration, we have $\E{\overline{m}\mid T_1\M}=\E{\overline{m}\mid T_2\M}=0$.
That is $(U_{\{1,2\}}^i\overline{m})_{i\in\Z^2}$ is an OMD random field. Moreover, using
(\ref{decomposition_m1barre_en_m2barre_g2barre}), we have
$$
\E{\overline{m}_2\mid T_1T_3\M}=\E{\overline{m}_1-(I-U_2)\overline{g}_2\mid T_1T_3\M}=-(I-
U_2)\E{\overline{g}_2\mid T_1T_2T_3\M}.
$$
Since $\overline{g}_2$ is $T_2T_3\M$-measurable and $(T^{-i}\mathcal{M})_{i\in\Z^d}$ is a commuting filtration, we have also
$\E{\overline{g}_2\mid T_1T_2T_3\M}=\E{\overline{g}_2\mid T_1\M}$ and consequently
\begin{equation}\label{esperance_m2barre_sachantT1T3M}
\E{\overline{m}_2\mid T_1T_3\M}=-(I-U_2)\E{\overline{g}_2\mid T_1\M}.
\end{equation}
Combining (\ref{decomposition_g3_en_m1barre_g1barre}),
(\ref{decomposition_m1barre_en_m2barre_g2barre}),
(\ref{definition_mbarre}) and (\ref{esperance_m2barre_sachantT1T3M}), we obtain
$$
g_3=\overline{m}+(I-U_1)\overline{g}_1+(I-U_2)[\overline{g}_2-\E{\overline{g}_2\mid T_1\M}].
$$
As before, since $\overline{g}_1$ is $T_1T_3\M$-measurable and $\sum_{k\geqslant
1}\norm{\E{\overline{g}_1\mid T_2^k\mathcal M}}_p<\infty$ (by Lemma \ref{lemme_cle}), we have
$$
\overline{g}_1=\overline{\overline{m}}_2+(I-U_2)\overline{\overline{g}}_2
$$
where $\overline{\overline{m}}_2\in \mathbb{L}^p(\Omega, T_1T_3\mathcal M,\mu)\ominus\mathbb{L}^p(\Omega, T_1T_2T_3\mathcal M,\mu)$ and
$\overline{\overline{g}}_2\in \mathbb{L}^p(\Omega, T_1T_2T_3\mathcal M,\mu)$. In particular,
$(U_2^k\overline{\overline{m}}_2)_{k\in\Z}$ is a MD sequence. Consequently,
\begin{align*}
f=m&+(I-U_1)[g_1-\E{g_1\mid T_2\M}-\E{g_1\mid T_3\M} +\E{g_1\mid T_2T_3\M}]\\
&+(I-U_2)[g_2-\E{g_2\mid T_3\M}]\\
&+(I-U_3)\overline{m}+(I-U_3)(I-U_2)[\overline{g}_2-\E{\overline{g}_2\mid T_1\M}]\\
&+(I-U_3)(I-U_1)\overline{\overline{m}}_2+(I-U_3)(I-U_1)(I-U_2)\overline{\overline{g}}_2.
\end{align*}
Since $g_2$ is $T_2\M$-measurable and $\sum_{k\geqslant 1}\norm{\E{g_2\mid T_1^k\mathcal
M}}_p<\infty$ (by Lemma \ref{lemme_cle}), we have
\begin{equation}\label{decomposition_g2_en_m1barrebarre_g1barrebarre}
g_2=\overline{\overline{m}}_1+(I-U_1)\overline{\overline{g}}_1
\end{equation}
where $\overline{\overline{m}}_1\in \mathbb{L}^p(\Omega, T_2\M,\mu)\ominus\mathbb{L}^p(\Omega,T_1 T_2\M,\mu)$ and $\overline{\overline{g}}_1\in
\mathbb{L}^p(\Omega, T_1T_2\M,\mu)$. Denoting
\begin{equation}\label{definition_mbarrebarre}
\overline{\overline{m}}:=\overline{\overline{m}}_1-\E{\overline{\overline{m}}_1\mid T_3\M},
\end{equation}
and applying Lemma \ref{lemme_cle}, we have $\E{\overline{\overline{m}}\mid
T_1\M}=\E{\overline{\overline{m}}\mid T_3\M}=0$. So, $(U_{\{1,3\}}^i\overline{\overline{m}})_{i\in\Z^2}$
is an OMD random field. Moreover, using (\ref{decomposition_g2_en_m1barrebarre_g1barrebarre}), we derive
\begin{equation}\label{esperance_m1barrebarre_sachant_T3M}
\E{\overline{\overline{m}}_1\mid T_3\M}=\E{g_2\mid T_3\M}-(I-U_1)\E{\overline{\overline{g}}_1\mid
T_1T_3\M}.
\end{equation}
Since $\overline{\overline{g}}_1$ is $T_1\M$-measurable and $(T^{-i}\mathcal{M})_{i\in\Z^3}$ is a commuting filtration,  we know that
$\E{\overline{\overline{g}}_1\mid T_1T_3\M}=\E{\overline{\overline{g}}_1\mid T_3\M}$. Combining
(\ref{decomposition_g2_en_m1barrebarre_g1barrebarre}), (\ref{definition_mbarrebarre}) and
(\ref{esperance_m1barrebarre_sachant_T3M}), we obtain
$$
g_2-\E{g_2\mid T_3\M}=\overline{\overline{m}}+(I-U_1)
[\overline{\overline{g}}_1-\E{\overline{\overline{g}}_1\mid T_3\M}].
$$
Finally,
\begin{align*}
f=m&+(I-U_1)[g_1-\E{g_1\mid T_2\M}-\E{g_1\mid T_3\M} +\E{g_1\mid T_2T_3\M}]\\
&+(I-U_2)\overline{\overline{m}}+(I-U_2)(I-U_1)
[\overline{\overline{g}}_1-\E{\overline{\overline{g}}_1\mid T_3\M}]\\
&+(I-U_3)\overline{m}+(I-U_3)(I-U_2)[\overline{g}_2-\E{\overline{g}_2\mid T_1\M}]\\
&+(I-U_3)(I-U_1)\overline{\overline{m}}_2+(I-U_3)(I-U_1)(I-U_2)\overline{\overline{g}}_2.
\end{align*}
The proof of Theorem \ref{decomposition_dim_d} is complete for $d=3$. The proof of Theorem
\ref{decomposition_dim_d} for $d\geqslant 4$ can be done in the same way. It is left to the reader.\\
\\
%
%
%
%
%
%
%
{\em Proof of Proposition \ref{Moment_inequalities}}. Let $(X_i)_{i\in\Z^d}$ be an OMD random field with respect to a commuting filtration $(\mathcal{F}_i)_{i\in\Z^d}$. Again, for simplicity, we consider only the case $d=2$. Let $n=(n_1,n_2)\succeq 0$ be fixed in $\Z^2$ and consider $(Y_i)_{i\in\Z}$ defined for any $i$ in $\Z$ by $Y_i=\sum_{j=0}^{n_2}X_{i,j}$. One can notice that $(Y_i)_{i\in\Z}$ is a MD sequence with respect to the filtration $(\vee_{j\in\Z}\mathcal{F}_{(i,j)})_{i\in\Z}$. Consequently, by Burkholder's inequality (see \cite{Hall--Heyde1980}, Theorem 2.10), we have
$$
\norm{\sum_{i=0}^{n_1}\sum_{j=0}^{n_2}X_{i,j}}_p
\leqslant\kappa\sqrt{p}\left(\sum_{i=0}^{n_1}
\norm{Y_i}_p^2\right)^{1/2}.
$$
Moreover, since for any $i$ in $\Z$, $(X_{i,j})_{j\in\Z}$ is a MD sequence with respect to the filtration $(\vee_{i\in\Z}\mathcal{F}_{(i,j)})_{j\in\Z}$, we have also
$$
\norm{Y_i}_p=\norm{\sum_{j=0}^{n_2}X_{i,j}}_p
\leqslant\kappa\sqrt{p}\left(\sum_{j=0}^{n_2}
\norm{X_{i,j}}_p^2\right)^{1/2}.
$$
Consequently, we obtain
\begin{equation}\label{Moment_inequality_dim2}
\norm{\sum_{i=0}^{n_1}\sum_{j=0}^{n_2}X_{i,j}}_p
\leqslant\kappa p\left(\sum_{i=0}^{n_1}\sum_{j=0}^{n_2}
\norm{X_{i,j}}_p^2\right)^{1/2}.
\end{equation}
In order to prove the optimality of the constant $p$ in ($\ref{Moment_inequality_dim2}$), arguing as in Wang and Woodroofe \cite{Wang--Woodroofe2011b} (Example 1, page 12), we consider a sequence $(\eta_i)_{i\in\Z}$ of iid real random variables satisfying $\mu(\eta_0=1)=\mu(\eta_0=-1)=1/2$. Let also $(\eta_{i}^{'})_{i\in\Z}$ be an independent copy of $(\eta_i)_{i\in\Z}$ and consider the filtrations $(\mathcal{G}_k)_{k\in\Z}$ and $(\mathcal{H}_k)_{k\in\Z}$ defined for any $k$ in $\Z$ by  $\mathcal{G}_k=\sigma(\eta_s;s\leqslant k)$ and $\mathcal{H}_k=\sigma(\eta_s^{'};s\leqslant k)$. For any $(i,j)$ in $\Z^2$, we denote  $Z_{i,j}=\eta_i\eta_{j}^{'}$. Then $(Z_{i,j})_{(i,j)\in\Z^2}$ is an OMD random field with respect to the commuting filtration $(\mathcal{F}_{i,j})_{(i,j)\in\Z^2}$ defined by $\mathcal{F}_{i,j}=\mathcal{G}_i\vee \mathcal{H}_j$ for any $(i,j)$ in $\Z^2$. Let $C$ be a positive constant such that for any $n=(n_1,n_2)\succeq 0$,
$$
\norm{\sum_{i=0}^{n_1}\eta_i}_p
\times\norm{\sum_{j=0}^{n_2}\eta_j^{'}}_p=
\norm{\sum_{i=0}^{n_1}\sum_{j=0}^{n_2}Z_{i,j}}_p
\leqslant C\left(\sum_{i=0}^{n_1}\sum_{j=0}^{n_2}
\norm{Z_{i,j}}_p^2\right)^{1/2}\leqslant C\sqrt{n_1n_2}.
$$
Applying the CLT for iid real random variables,
we derive $C\geqslant \norm{N}_p^2$ where $N$ is a standard normal random variable. Since there exists $\kappa>0$ such that $\norm{N}_p^2\geqslant\kappa p$, we derive (\ref{borne_inf_constants}). The proof of Proposition \ref{Moment_inequalities} is complete.
%
%
%
%
%
%
\\
\\
{\em Proof of Theorem \ref{IP_OMD}}. The convergence of the finite-dimensional laws of the partial sums process is a direct consequence of the CLT by Wang and Woodroofe (\cite{Wang--Woodroofe2011b}, Theorem 5). So, it suffices to establish the tightness property of the partial sums process. Assume that $\mathcal{A}$ is a VC-class with index $V$ and there exists $p>2(V-1)$ such that $X_0$ belongs to $\L^p(\Omega,\F,\mu)$ and \eqref{criterion_IP_dim_d_norm_p} holds. Then there exists
a positive constant $K$ such that for any $0<\varepsilon<1$,
we have (see Van der Vaart and Wellner \cite{van-der-Vaart-Wellner}, Theorem 2.6.4)
$$
N(\A,\rho,\varepsilon)\leqslant KV(4e)^V\left(\frac{1}{\varepsilon}\right)^{2(V-1)}
$$
where $N(\A,\rho,\varepsilon)$ is the smallest number of open balls of radius $\epsilon$ with respect to $\rho$
which form a covering of $\A$. Since $p>2(V-1)$, we have
\begin{equation}\label{metric_entropy_Lp}
\int_0^1\left(N(\A,\rho,\varepsilon)\right)^{\frac{1}{p}}d\varepsilon <\infty.
\end{equation}
Moreover, using Proposition $\ref{Moment_inequality_under_proj_condition}$, we derive
\begin{equation}\label{IP_Lipschitz_Lp}
\norm{n^{-d/2}S_n(A)-n^{-d/2}S_n(B)}_p\leqslant\kappa p^{d/2}\rho(A,B)\norm{X_0}_p
\end{equation}
for any positive integer $n$ and any $A$ and $B$ in $\A$.
Combining ($\ref{metric_entropy_Lp}$) and \eqref{IP_Lipschitz_Lp} and
applying Theorem 11.6 in Ledoux and Talagrand \cite{Led-Tal}, we obtain that for each positive $\epsilon$ there exists a positive real $\delta$, depending on $\epsilon$ and on the value of
the entropy integral ($\ref{metric_entropy_Lp}$) but not on $n$, such that
\begin{equation}\label{esperance_plus_petite_epsilon}
\E{\sup_{\stackrel{A,B\in\A}{\rho(A,B)<\delta}}\vert n^{-d/2}S_{n}(A)-n^{-d/2}S_{n}(B)\vert}<\epsilon.
\end{equation}
In particular, for any $x>0$, we have
$$
\lim_{\delta\to0}\limsup_{n\to\infty}
\mu\left(\sup_{\substack{A,B\in\A \\
\rho(A,B)<\delta}}\big\vert
n^{-d/2}S_{n}(A)-n^{-d/2}S_{n}(B)\big\vert>x\right)=0.
$$
Consequently, the partial sum process $\{n^{-d/2}S_n(A);A\in\A\}_{n\geqslant 1}$ is tight in the space $\mathcal{C(\A)}$ and the WIP holds. The proof of Theorem \ref{IP_OMD} is complete.\\
\\
{\em Proof of Theorem \ref{IP_OMD_bis}}.  Let $(\Omega, \F, \mu, \{T^k\}_{k\in\Z^d})$ be a dynamical system
(i.e. $(\Omega, \F,\mu)$ is a probability space and $T^k:\Omega\to\Omega$ is a measure-preserving
transformation for any $k$ in $\Z^d$ satisfying $T^i\circ T^j=T^{i+j}$ for any $i$ and any $j$ in
$\Z^d$) and let  $(\varepsilon_i)_{i\in\Z^d}$ be a field of iid real random variables defined on
$(\Omega, \F,\mu)$. Let $\M\subset \F$ be the $\sigma$-algebra generated by the random variables
$\varepsilon_i$ for $i\preceq 0$ and let $f:\Omega\to\R$ be $\mathcal M$-measurable. We consider the
stationary real random field $(f\circ T^i)_{i\in\Z^d}$ and the partial sum process $\left\{S_n(f,t)\,;
\,t\in[0,1]^d\right\}_{n\geqslant 1}$ defined for any integer $n\geqslant 1$ and any $t$ in $[0,1]^d$ by
\begin{equation}\label{def_Snft}
S_n(f,t)=\sum_{i\in \langle n\rangle^d}\lambda([0,nt]\cap R_i)f\circ T^i
\end{equation}
where $\lambda$ is the Lebesgue measure on $\R^d$ and $R_i=]i_1-1,i_1]\times ...\times]i_d-1,i_d]$ is
the unit cube with upper corner $i=(i_1,...,i_d)$ in $\langle n\rangle^d$. Again, the convergence of the finite-dimensional laws of the process $\left\{n^{-d/2}S_n(f,t)\,;
\,t\in[0,1]^d\right\}_{n\geqslant 1}$ is a
direct consequence of the CLT established by Wang and Woodroofe (\cite{Wang--Woodroofe2011b}, Theorem
3.2). In order to obtain the tightness property of the partial sum process, it suffices to
establish for any $\varepsilon>0$,
$$
\lim_{\delta\to 0}\limsup_{n\to\infty}
\mu\left(\sup_{\substack{s,t\in[0,1]^d \\ \vert s-t\vert<\delta}}n^{-d/2}\vert S_n(f,s)-
S_n(f,t)\vert>\varepsilon\right)=0
$$
where $\vert x\vert=\max_{k\in\langle d\rangle}\vert x_k\vert$ for any $x=(x_1,...,x_d)$ in $[0,1]^d$. For
simplicity,
we are going to consider only the case $d=2$. By Theorem \ref{decomposition_dim_d}, we have
\begin{equation}
f=m+(I-U_1)m_1+(I-U_2)m_2+(I-U_1)(I-U_2)g,
\end{equation}
where $m$, $m_1$, $m_2$ and $g$ are square-integrable functions defined on $\Omega$ such that
$(U^im)_{i\in\Z^2}$ is an OMD random field and $(U_2^k m_1)_{k\in\Z}$ and $(U_1^k m_2)_{k\in\Z}$ are MD
sequences. In the sequel, for any real $x$, we denote by $[x]$ the integer part of $x$. Let $n\geqslant 1$
and $t=(t_1,t_2)$ in $[0,1]^2$. For any $1\leqslant i\leqslant [nt_1]+1$ and any $1\leqslant j\leqslant [nt_2]+1$, we denote
$\lambda_{i,j}(t)=\lambda\left([0,nt]\cap R_{(i,j)}\right)$. We have
$$
S_n((I-U_1)m_1,t)
=\sum_{i=1}^{[nt_1]+1}\sum_{j=1}^{[nt_2]+1}\lambda_{i,j}(t)U^{(i,j)}(I-U_1)m_1
=\sum_{j=1}^{[nt_2]+1}U_2^j\sum_{i=1}^{[nt_1]+1}\lambda_{i,j}(t)\left(U_1^im_1-U_1^{i+1}m_1\right).
$$
Using Abel's transformation and noting that $\lambda_{i+1,j}(t)=\lambda_{i,j}(t)$ for any $1\leqslant i\leqslant
[nt_1]-1$ and any $1\leqslant j\leqslant [nt_2]+1$, we obtain that $S_n((I-U_1)m_1,t)$ equals
\begin{align*}
&\sum_{j=1}^{[nt_2]+1}U_2^j\left\{\lambda_{[nt_1]+1,j}(t)\left(U_1m_1-
U_1^{[nt_1]+2}m_1\right)-\sum_{i=1}^{[nt_1]}\left(U_1m_1-U_1^{i+1}m_1\right)\left(\lambda_{i+1,j}
(t)-\lambda_{i,j}(t)\right)\right\}\\
&=\sum_{j=1}^{[nt_2]+1}U_2^j\left\{\lambda_{[nt_1]+1,j}(t)\left(U_1m_1-
U_1^{[nt_1]+2}m_1\right)-\left(U_1m_1-U_1^{[nt_1]+1}m_1\right)\left(\lambda_{[nt_1]+1,j}
(t)-\lambda_{[nt_1],j}(t)\right)\right\}\\
&=U_1(I-U_1^{[nt_1]+1})\sum_{j=1}^{[nt_2]+1}\lambda_{[nt_1]+1,j}(t)\,U_2^jm_1
-U_1(I-U_1^{[nt_1]})\sum_{j=1}^{[nt_2]+1}\left(\lambda_{[nt_1]+1,j}(t)
-\lambda_{[nt_1],j}(t)\right)U_2^jm_1.
\end{align*}
Moreover, since $\lambda_{i,j}(t)=\lambda_{i,1}(t)$ for any $1\leqslant i\leqslant [nt_1]+1$ and any $1\leqslant j\leqslant
[nt_2]$, we derive
\begin{align*}
S_n((I-U_1)m_1,t)&=U_1(I-U_1^{[nt_1]+1})\lambda_{[nt_1]+1,1}(t)\sum_{j=1}^{[nt_2]}U_2^jm_1\\
&\quad+U_1(I-U_1^{[nt_1]+1})\lambda_{[nt_1]+1,[nt_2]+1}(t)\,U_2^{[nt_2]+1}m_1\\
&\quad-U_1(I-U_1^{[nt_1]})\left(\lambda_{[nt_1]+1,1}(t)
-\lambda_{[nt_1],1}(t)\right)\sum_{j=1}^{[nt_2]}U_2^jm_1\\
&\quad-U_1(I-U_1^{[nt_1]})\left(\lambda_{[nt_1]+1,[nt_2]+1}(t)
-\lambda_{[nt_1],[nt_2]+1}(t)\right)U_2^{[nt_2]+1}m_1.
\end{align*}
So, we obtain
\begin{equation}\label{inequality_sup}
\sup_{t\in [0,1]^2}\vert S_n((I-U_1)m_1,t)\vert
\leqslant 4\max_{1\leqslant l,k\leqslant n+2}U_1^{l}U_2^k\vert m_1\vert+4\max_{1\leqslant l,k\leqslant
n+2}U_1^{l}\left\vert\sum_{j=1}^{k}U_2^jm_1\right\vert.
\end{equation}
Let $x>0$ be fixed. Since $m_1\in\L^2(\Omega,\F,\mu)$, we have
\begin{equation}\label{max1}
\mu\left(\max_{1\leqslant l,k\leqslant n+2}U_1^{l}U_2^k\vert m_1\vert>nx\right)\leqslant \kappa
n^2\mu\left(m_1^2>n^2x^2\right)\converge{n}{\infty}{ }0.
\end{equation}
In the other part,
\begin{equation}\label{max2}
\mu\left(\max_{1\leqslant l,k\leqslant n+2}U_1^{l}\left\vert\sum_{j=1}^{k}U_2^jm_1\right\vert>xn\right)
=\mu\left(\max_{1\leqslant l\leqslant n+2}U_1^l\left(\frac{1}{\sqrt{n}}\max_{1\leqslant k\leqslant
n+2}\sum_{j=1}^{k}U_2^jm_1\right)^2>nx^2\right).
\end{equation}

\begin{Lem}\label{ui}
Let $(Z_n)_{n\geqslant 1}$ be a sequence of uniformly integrable real random variables. For any $s$ in
$\langle d\rangle$,
$$
\limsup_{n\to\infty}\mu\left(\max_{1\leqslant i_1,..,i_s\leqslant n}U_1^{i_1}\dots
U_s^{i_s}\abs{Z_n}>n^s\right)=0.
$$
\end{Lem}
{\em Proof of Lemma \ref{ui}}. Let $n$ be a positive integer. For any $s$ in $\langle d\rangle$, we
denote
$$
p_n(s):=\mu\left(\max_{1\leqslant i_1,\dots,i_s\leqslant n}U_1^{i_1}\dots U_s^{i_s}\abs{Z_n}>n^s\right).
$$
Let $R$ be a positive real number. We have
$$
p_n(s)\leqslant \frac{2R}{n^s}+n^s\,\mu\left(\abs{Z_n}\ind{\{\abs{Z_n}> R\}}>\frac{n^s}{2}\right)
\leqslant\frac{2R}{n^s}+2\sup_{k\geqslant 1}\E{\abs{Z_k}\ind{\{\abs{Z_k}> R\}}}.
$$
Consequently, $\limsup_{n\to\infty}p_n(s)\leqslant2\sup_{k\geqslant 1}\E{\abs{Z_k}\ind{\{\abs{Z_k}> R\}}}\converge{R}{\infty}{}0$. The
proof of Lemma \ref{ui} is complete.
\begin{Lem}\label{ui_bis}
The sequence $\left\{\left(\frac{1}{\sqrt{n}}\max_{1\leqslant k\leqslant n+2}\sum_{j=1}^{k}U_2^jm_1\right)^2;
\,n\geqslant 1\right\}$ is uniformly integrable.
\end{Lem}
{\em Proof of Lemma \ref{ui_bis}}. Since $(U_2^km_1)_{k\in\Z}$ is a MD sequence, using Doob's
inequality, we derive
$$
\norm{\max_{1\leqslant k\leqslant n+2}\sum_{j=1}^{k}U_2^jm_1}_2\leqslant 2
\norm{\sum_{j=1}^{n+2}U_2^jm_1}_2\leqslant \kappa\sqrt{n}\norm{m_1}_2.
$$
So, $\left\{\left(\frac{1}{\sqrt{n}}\max_{1\leqslant k\leqslant n+2}\sum_{j=1}^{k}U_2^jm_1\right)^2;\,n\geqslant
1\right\}$ is bounded in $\L^1(\Omega,\F,\mu)$. Let $M$ be a fixed positive constant. We have
$m_1=m_1^{'}+m_1^{''}$ where
\begin{align*}
m_1^{'}&=m_1\ind{\vert m_1\vert\leqslant M}-\E{m_1\ind{\vert m_1\vert\leqslant M}\mid T_2\mathcal{M}}\\
m_1^{''}&=m_1\ind{\vert m_1\vert> M}-\E{m_1\ind{\vert m_1\vert>M}\mid T_2\mathcal{M}}.
\end{align*}
Moreover, if $A$ belongs to $\F$ then
\begin{align*}
\int_A \left(\frac{1}{\sqrt{n}}\max_{1\leqslant k\leqslant n+2}\sum_{j=1}^{k}U_2^jm_1\right)^2d\mu
&\leqslant2\int_A\left(\frac{1}{\sqrt{n}}\max_{1\leqslant k\leqslant
n+2}\sum_{j=1}^{k}U_2^jm_1^{'}\right)^2d\mu\\
&\quad+2\int_A\left(\frac{1}{\sqrt{n}}\max_{1\leqslant k\leqslant
n+2}\sum_{j=1}^{k}U_2^jm_1^{''}\right)^2d\mu.
\end{align*}
Since $(U_2^km_1^{'})_{k\in\Z}$ and $(U_2^km_1^{''})_{k\in\Z}$ are MD sequences, using Schwarz's
inequality, we obtain
\begin{align*}
\int_A\left(\frac{1}{\sqrt{n}}\max_{1\leqslant k\leqslant n+2}\sum_{j=1}^{k}U_2^jm_1\right)^2d\mu
&\leqslant2\norm{\frac{1}{\sqrt{n}}\max_{1\leqslant k\leqslant n+2}\sum_{j=1}^{k}U_2^jm_1^{'}}_4^2\sqrt{\mu(A)}\\
&\quad+2\norm{\frac{1}{\sqrt{n}}\max_{1\leqslant k\leqslant n+2}\sum_{j=1}^{k}U_2^jm_1^{''}}_2^2.
\end{align*}
Keeping in mind that $m_1^{'}$ is bounded by $M$ and using again Doob's inequality, there exists a
positive constant $\kappa_0$ such that
$$
\int_A\left(\frac{1}{\sqrt{n}}\max_{1\leqslant k\leqslant n+2}\sum_{j=1}^{k}U_2^jm_1\right)^2d\mu
\leqslant\kappa_0 \left(M^2\sqrt{\P(A)}+\E{m_1^2\ind{\vert m_1\vert>M}}\right).
$$
Let $\varepsilon>0$ be fixed and let $M>0$ such that $\kappa_0\E{m_1^2\ind{\vert
m_1\vert>M}}\leqslant\frac{\varepsilon}{2}$. One can choose the measurable set $A$ in $\F$ such that
$\kappa_0 M^2\sqrt{\mu(A)}\leqslant\frac{\varepsilon}{2}$ and consequently
$$
\sup_{n\geqslant 1}\int_A\left(\frac{1}{\sqrt{n}}\max_{1\leqslant k\leqslant
n+2}\sum_{j=1}^{k}U_2^jm_1\right)^2d\mu\leqslant\varepsilon.
$$
The proof of Lemma \ref{ui_bis} is complete.\\
\\
Combining (\ref{inequality_sup}), (\ref{max1}), (\ref{max2}),  Lemma \ref{ui} and Lemma \ref{ui_bis},
we obtain
\begin{equation}\label{limsup_U1}
\limsup_{n\to\infty}\mu\left(\sup_{t\in[0,1]^2}\vert S_n((I-U_1)m_1,t)\vert>xn\right)=0.
\end{equation}
In a similar way, we derive also
\begin{equation}\label{limsup_U2}
\limsup_{n\to\infty}\mu\left(\sup_{t\in[0,1]^2}\vert S_n((I-U_2)m_1,t)\vert>xn\right)=0.
\end{equation}
Now, noting that $\lambda_{i,j}(t)=\lambda_{i,1}(t)$ for any $1\leqslant i\leqslant [nt_1]+1$ and any $1\leqslant
j\leqslant [nt_2]$, we have $S_n((I-U_1)(I-U_2)m_1,t)$ equals
\begin{align*}
&\sum_{i=1}^{[nt_1]+1}\sum_{j=1}^{[nt_2]+1}
\lambda_{i,j}(t)\,U^{(i,j)}(I-U_1)(I-U_2)m_1\\
&=\sum_{i=1}^{[nt_1]+1}U_1^i(I-U_1)
\left(\lambda_{i,1}(t)\sum_{j=1}^{[nt_2]}(U_2^j-U_2^{j+1})m_1+\lambda_{i,[nt_2]+1}(t)U_2^{[nt_2]+1}(I-
U_2)m_1\right)\\
&=U_2(I-U_2^{[nt_2]})\sum_{i=1}^{[nt_1]+1}\lambda_{i,1}(t)(U_1^i-U_1^{i+1})m_1+U_2^{[nt_2]+1}(I-
U_2)\sum_{i=1}^{[nt_1]+1}\lambda_{i,[nt_2]+1}(t)(U_1^im_1-U_1^{i+1}m_1).
\end{align*}
Since $\lambda_{i,j}(t)=\lambda_{1,j}(t)$ for any $1\leqslant i\leqslant [nt_1]$ and any $1\leqslant j\leqslant [nt_2]+1$,
we derive
\begin{align*}
S_n((I-U_1)(I-U_2)m_1,t)&=\lambda_{1,1}(t)U_2(I-U_2^{[nt_2]})U_1(I-U_1^{[nt_1]})m_1\\
&\quad+\lambda_{[nt_1]+1,1}(t)U_2(I-U_2^{[nt_2]})U_1^{[nt_1]+1}\left(I-U_1\right)m_1\\
&\quad+\lambda_{1,[nt_2]+1}(t)U_2^{[nt_2]+1}(I-U_2)U_1(I-U_1^{[nt_1]+1})m_1\\
&\quad+\lambda_{[nt_1]+1,[nt_2]+1}(t)U_2^{[nt_2]+1}(I-U_2)U_1^{[nt_1]+1}(I-U_1)m_1.
\end{align*}
Thus
$$
\sup_{t\in[0,1]^2}\vert S_n((I-U_1)(I-U_2)m_1,t)\vert
\leqslant\kappa\max_{1\leqslant k,l\leqslant n+2}U_1^kU_2^l\vert m_1\vert$$
and for any positive $x$,
\begin{equation}\label{limsup_cobord}
\mu\left(\sup_{t\in[0,1]^2}\vert S_n((I-U_1)(I-U_2)m_1,t)\vert>xn\right)\leqslant\kappa
n^2\mu\left(m_1^2>n^2x^2\right)\converge{n}{\infty}{}0.
\end{equation}
Now, it suffices to prove the tightness of the process $\{\frac{1}{n} S_n(m,t);t\in[0,1]^2\}_{n\geqslant 1}$. That is,
for any positive $x$,
\begin{equation}\label{tension_orthomartingale}
\lim_{\delta\to 0}\limsup_{n\to\infty}\mu\left(
\sup_{\substack{s,t\in[0,1]^2\\ \vert s-t\vert<\delta}}\vert S_n(m,s)-S_n(m,t)\vert>xn\right)=0.
\end{equation}
Let $n$ be a positive integer and let $s=(s_1,s_2)$ and $t=(t_1,t_2)$ be fixed in $[0,1]^2$. We denote
$\Delta_n(s,t)=S_n(m,s)-S_n(m,t)$ and for any $i$ and $j$ in $\langle n\rangle$,
$$
\beta_{i,j}=\lambda_{i,j}(s)-\lambda_{i,j}(t)=\lambda\left([0,ns]\cap R_{(i,j)}\right)-\lambda\left([0,nt]\cap R_{(i,j)}\right).
$$
Noting that $\beta_{i,j}=0$ for any $1\leqslant i\leqslant [n(s_1\wedge t_1)]$
and any $1\leqslant j\leqslant [n(s_1\wedge
t_1)]$, we have $
\Delta_n(s,t)=\Delta_n^{'}(s,t)+\Delta_n^{''}(s,t)$ where
$$
\Delta_n^{'}(s,t)=\sum_{i=[n(s_1\wedge t_1)]+1}^{[n(s_1\vee t_1)]+1}\sum_{j=1}^{[n(s_2\wedge
t_2)]+1}\beta_{i,j}\,U^{(i,j)}m
\quad\textrm{and}\quad
\Delta_n^{''}(s,t)=\sum_{i=1}^{[n(s_1\wedge t_1)]+1}\sum_{j=[n(s_2\wedge t_2)]+1}^{[n(s_2\vee
t_2)]+1}\beta_{i,j}\,U^{(i,j)}m.
$$
Moreover, $\Delta_n^{'}(s,t)=\Delta_{1,n}^{'}(s,t)+\Delta_{2,n}^{'}(s,t)+\Delta_{3,n}^{'}
(s,t)+\Delta_{4,n}^{'}(s,t)$ where
\begin{align*}
\Delta_{1,n}^{'}(s,t)&=\sum_{i=[n(s_1\wedge t_1)]+2}^{[n(s_1\vee t_1)]}\sum_{j=1}^{[n(s_2\wedge
t_2)]}\beta_{i,j}\,U^{(i,j)}m\\
\Delta_{2,n}^{'}(s,t)&=\sum_{j=1}^{[n(s_2\wedge t_2)]}\beta_{[n(s_1\vee t_1)]+1,j}\,U^{([n(s_1\vee
t_1)]+1,j)}m\\
\Delta_{3,n}^{'}(s,t)&=\sum_{j=1}^{[n(s_2\wedge t_2)]}\beta_{[n(s_1\wedge t_1)]+1,j}\,U^{([n(s_1\wedge
t_1)]+1,j)}m\\
\Delta_{4,n}^{'}(s,t)&=\sum_{i=[n(s_1\wedge t_1)]+1}^{[n(s_1\vee t_1)]+1}\beta_{i,[n(s_2\wedge
t_2)]+1}\,U^{(i,[n(s_2\wedge t_2)]+1)}m.
\end{align*}
Let $\alpha$ in $\{-1,+1\}$ such that $\beta_{i,j}=\alpha$ if $[n(s_1\wedge t_1)]+2\leqslant i\leqslant
[n(s_1\vee t_1)]$ and $1\leqslant j\leqslant [n(s_2\wedge t_2)]$
So,
$$
\Delta_{1,n}^{'}(s,t)=\alpha\sum_{i=[n(s_1\wedge t_1)]+2}^{[n(s_1\vee t_1)]}\sum_{j=1}^{[n(s_2\wedge
t_2)]}U^{(i,j)}m
$$
and for any positive $x$,
\begin{align*}
\mu\left(\sup_{\substack{s,t\in[0,1]^2\\ \vert s-t\vert<\delta}}\vert\Delta_{1,n}^{'}
(s,t)\vert>nx\right)
&\leqslant \sum_{k=0}^{[\frac{1}{\delta}]}\mu\left(\max_{\substack{1\leqslant p\leqslant n\\ r\in[0,\delta]}}\left\vert
\sum_{i=[nk\delta]+2}^{[n(k\delta+r)]}\sum_{j=1}^pU^{(i,j)}m
\right\vert>nx\right)\\
&=\sum_{k=0}^{[\frac{1}{\delta}]}\mu\left(\max_{\substack{1\leqslant p\leqslant n\\ r\in[0,\delta]}}\left\vert
\sum_{i=1}^{[n(k\delta+r)]-[nk\delta]-1}\sum_{j=1}^pU^{(i,j)}m\right\vert>n
x\right).
\end{align*}
Since $[n(k\delta+r)]-[nk\delta]-1$ is an integer smaller than $[nr]$, we obtain
\begin{align*}
\mu\left(\sup_{\substack{s,t\in[0,1]^2\\ \vert s-t\vert<\delta}}\vert\Delta_{1,n}^{'}
(s,t)\vert>nx\right)
&\leqslant\left(1+\frac{1}{\delta}\right)\mu\left(\max_{\substack{1\leqslant p\leqslant n\\ 1\leqslant q\leqslant
[n\delta]}}\left\vert\sum_{i=1}^q
\sum_{j=1}^pU^{(i,j)}m\right\vert>nx\right)\\
&=\left(1+\frac{1}{\delta}\right)\mu\left(\max_{\substack{1\leqslant p\leqslant n\\ 1\leqslant q\leqslant
[n\delta]}}\left(\frac{1}{n\sqrt{\delta}}\sum_{i=1}^q
\sum_{j=1}^pU^{(i,j)}m\right)^2>\frac{x^2}{\delta}\right)\\
&\leqslant \left(\frac{1+\delta}{x^2}\right)\mathbb{E}_{\frac{x^2}
{\delta}}\left(\max_{\substack{1\leqslant p\leqslant n\\ 1\leqslant q\leqslant [n\delta]}}\left(\frac{1}
{n\sqrt{\delta}}\sum_{i=1}^q
\sum_{j=1}^pU^{(i,j)}m\right)^2\right)
\end{align*}
where we used the notation $\mathbb{E}_{A}(Z)=\E{Z\ind{\vert Z\vert>A}}$ for any $A>0$ and any $Z$ in $\L^1(\Omega,\F,\mu)$.
\begin{Lem}\label{ui_orthomartingale}
The family $\left\{\max_{\substack{1\leqslant p\leqslant n\\ 1\leqslant q\leqslant [n\delta]}}\left(\frac{1}
{n\sqrt{\delta}}\sum_{i=1}^q
\sum_{j=1}^pU^{(i,j)}m\right)^2\,;\,n\geqslant 1, \delta>0\right\}$ is uniformly integrable.
\end{Lem}
{\em Proof of Lemma \ref{ui_orthomartingale}}. The proof follows the same lines as the proof of
Lemma \ref{ui_bis} using Cairoli's maximal inequality for orthomartingales (see
\cite{Khoshnevisan2002},
Theorem 2.3.1) instead of Doob's inequality for martingales. The proof of Lemma
\ref{ui_orthomartingale}
is complete.\\
\\
So, we obtain
\begin{equation}\label{tension_delta_1n_prime}
\lim_{\delta\to0}\limsup_{n\to\infty}
\mu\left(\sup_{\substack{s,t\in[0,1]^2\\ \vert s-t\vert<\delta}}\vert\Delta_{1,n}^{'}
(s,t)\vert>nx\right)=0.
\end{equation}
In the other part, since $\beta_{[n(s_1\vee t_1)]+1,j}=\beta_{[n(s_1\vee t_1)]+1,1}$ for any $1\leqslant
j\leqslant[n(s_2\wedge t_2)]$, we have
$$
\Delta_{2,n}^{'}(s,t)=\beta_{[n(s_1\vee t_1)]+1,1}U_1^{[n(s_1\vee t_1)]+1}\sum_{j=1}^{[n(s_2\wedge
t_2)]}U_2^jm
$$
and consequently
$$
\sup_{\substack{s,t\in[0,1]^2\\ \vert s-t\vert<\delta}}
\vert\Delta_{2,n}^{'}(s,t)\vert\leqslant\max_{\substack{1\leqslant k\leqslant n+1\\ 1\leqslant l\leqslant
n}}U_1^k\left\vert\sum_{j=1}^lU_2^jm\right\vert.
$$
So,
\begin{equation}\label{Inequality_Delta_1n_prime}
\mu\left(\sup_{\substack{s,t\in[0,1]^2\\ \vert s-t\vert<\delta}}
\vert\Delta_{2,n}^{'}(s,t)\vert>nx\right)
\leqslant \mu\left(\max_{1\leqslant k\leqslant n+1}U_1^k\left(\max_{1\leqslant l\leqslant n}\frac{1}
{\sqrt{n}}\left\vert\sum_{j=1}^lU_2^jm\right\vert
\right)^2>nx^2\right).
\end{equation}
Since $(U_2^km)_{k\in\Z}$ is a MD sequence, arguing as in Lemma \ref{ui_bis}, the
sequence  $\left\{\left(\max_{1\leqslant l\leqslant n}\frac{1}{\sqrt{n}}\left\vert\sum_{j=1}^lU_2^jm\right\vert
\right)^2\right\}_{n\geqslant 1}$ is uniformly integrable. Combining (\ref{Inequality_Delta_1n_prime}) and
Lemma \ref{ui}, we derive that for any $\delta>0$,
\begin{equation}\label{tension_delta_2n_prime}
\limsup_{n\to\infty}\mu\left(\sup_{\substack{s,t\in[0,1]^2\\ \vert s-t\vert<\delta}}
\vert\Delta_{2,n}^{'}(s,t)\vert>nx\right)=0.
\end{equation}
Similarly, we have also
\begin{equation}\label{tension_delta_3n_prime}
\limsup_{n\to\infty}\mu\left(\sup_{\substack{s,t\in[0,1]^2\\ \vert s-t\vert<\delta}}
\vert\Delta_{3,n}^{'}(s,t)\vert>nx\right)=0
\end{equation}
for any $\delta>0$. Moreover, for any $[n(s_1\wedge t_1)]+1\leqslant i\leqslant [n(s_1\vee t_1)]$, we have $\beta_{i,[n(s_2\wedge t_2)]+1}=\beta_{[n(s_1\wedge t_1)]+1,[n(s_2\wedge t_2)]+1}$ and consequently
\begin{align*}
\Delta_{4,n}^{'}(s,t)&=\beta_{[n(s_1\wedge t_1)]+1,[n(s_2\wedge t_2)]+1}\,U_2^{[n(s_2\wedge
t_2)]+1}\sum_{i=[n(s_1\wedge t_1)]+1}^{[n(s_1\vee t_1)]}U_1^{i}m\\
&\quad+\beta_{[n(s_1\vee t_1)]+1,[n(s_2\wedge t_2)]+1}\,U^{([n(s_1\vee t_1)]+1,[n(s_2\wedge t_2)]+1)}m
\end{align*}
and
\begin{align*}
\mu\left(\sup_{\substack{s,t\in[0,1]^2\\ \vert s-t\vert<\delta}}
\vert\Delta_{4,n}^{'}(s,t)\vert>nx\right)
&\leqslant
\mu\left(\max_{1\leqslant k\leqslant n+1}U_2^k\left(\max_{1\leqslant l\leqslant [n\delta]}\frac{1}
{\sqrt{n\delta}}\left\vert\sum_{j=1}^lU_1^jm\right\vert
\right)^2>\frac{nx^2}{2\delta}\right)\\
&\quad+2n^2\mu\left( m^2>\frac{n^2x^2}{4}\right)
\end{align*}
Arguing as in Lemma \ref{ui}, the family $\left\{\left(\max_{1\leqslant l\leqslant [n\delta]}\frac{1}
{\sqrt{n\delta}}\left\vert\sum_{j=1}^lU_1^jm\right\vert
\right)^2\,;\,n\geqslant 1, \delta>0\right\}$ is uniformly integrable since $(U_1^km)_{k\in\Z}$ is a MD sequence. By Lemma \ref{ui_bis}, we obtain for any $\delta>0$,
$$
\limsup_{n\to\infty}\mu\left(\max_{1\leqslant k\leqslant n+1}U_2^k\left(\max_{1\leqslant l\leqslant [n\delta]}\frac{1}
{\sqrt{n\delta}}\left\vert\sum_{j=1}^lU_1^jm\right\vert
\right)^2>\frac{nx^2}{2\delta}\right)=0.
$$
Moreover, $n^2\mu\left( m^2>\frac{n^2x^2}{4}\right)$ goes to zero as $n$ goes to infinity since $m$
belongs to $\L^2(\Omega,\F,\mu)$. Consequently, for any $\delta>0$,
\begin{equation}\label{tension_delta_4n_prime}
\limsup_{n\to\infty}\mu\left(\sup_{\substack{s,t\in[0,1]^2\\ \vert s-t\vert<\delta}}
\vert\Delta_{4,n}^{'}(s,t)\vert>nx\right)=0.
\end{equation}
Combining (\ref{tension_delta_1n_prime}), (\ref{tension_delta_2n_prime}),(\ref{tension_delta_3n_prime})
and (\ref{tension_delta_4n_prime}), we obtain
\begin{equation}\label{tension_delta_n_prime}
\lim_{\delta\to 0}\limsup_{n\to\infty}\mu\left(\sup_{\substack{s,t\in[0,1]^2\\ \vert s-t\vert<\delta}}
\vert\Delta_{n}^{'}(s,t)\vert>nx\right)=0.
\end{equation}
Similarly, one can check that
\begin{equation}\label{tension_delta_n_second}
\lim_{\delta\to 0}\limsup_{n\to\infty}\mu\left(\sup_{\substack{s,t\in[0,1]^2\\ \vert s-t\vert<\delta}}
\vert\Delta_{n}^{''}(s,t)\vert>nx\right)=0.
\end{equation}
Finally, keeping in mind $\Delta_n(s,t)=\Delta_n^{'}(s,t)+\Delta_n^{''}(s,t)$ and combining
(\ref{tension_delta_n_prime}) and (\ref{tension_delta_n_second}), we obtain
(\ref{tension_orthomartingale}). The proof of Theorem \ref{IP_OMD_bis} is complete.\\
\\
{\em Proof of Proposition \ref{IP_linear_random_fields}}. We shall use Rosenthal's inequality (\cite{Hall--Heyde1980}, Theorem~2.12). Let $p\geqslant 2$ be fixed. There exists a constant $C$ depending only on $p$ such that if $(Y_j)_{j\geqslant 1}$ is a
 sequence of independent zero-mean random variables and $n$ a positive integer then
 \begin{equation}\label{Rosenthal}
  \frac 1C\left(\sum_{j=1}^n\mathbb E[Y_j^2]\right)^{p/2}+\frac 1C\sum_{j=1}^n\mathbb E|Y_j|^p\leqslant
  \mathbb E\left|\sum_{j=1}^nY_j\right|^p\leqslant
  C\left(\sum_{j=1}^n\mathbb E[Y_j^2]\right)^{p/2}+C\sum_{j=1}^n\mathbb E|Y_j|^p.
 \end{equation}
Keeping in mind that $\Lambda_{k,s}=\{i=(i_1,..,i_d)\in\Z^d\,;\,i_s\geqslant k\}$ for any $k\geqslant 1$ and any $s$ in $\langle d\rangle$, we have $\E{X_0\mid\mathcal F_{k,s}}=\sum_{i\in\Lambda_{k,s}}a_i\varepsilon_{-i}$. Since $(\varepsilon_i)_{i\in\Z^d}$ is an iid real random field with $\varepsilon_0$ in $\L^p(\Omega,\F,\mu)$, we apply \eqref{Rosenthal} and the series $\sum_{k\geqslant 1}k^{d-1}
\lVert\mathbb E[X_0\mid\mathcal F_{k,s}]\rVert_p$ is convergent if and only if
$$
\sum_{k\geqslant 1}k^{d-1}\left\{
\left(\sum_{i\in \Lambda_{k,s}}|a_i|^2\right)^{1/2}+\left(\sum_{i\in \Lambda_{k,s}}|a_i|^p\right)^{1/p}\right\}<\infty.
$$
The result follows from the fact that $\sum_{i\in \Z}|c_i|^p\leqslant
\left(\sum_{i\in\Z}|c_i|^2\right)^{p/2}$ for any sequence $(c_i)_{i\in\Z}$ of real numbers. The proof of Proposition \ref{IP_linear_random_fields} is complete.\\
\\
{\em Proof of Theorem \ref{IP_Holder}}. We shall use Theorem~1 in \cite{Klicnarova2007} which states that
if a sequence of random processes $\{Y_n(t)\,;\,t\in [0,1]^d\}_{n\geqslant 1}$ whose
finite dimensional distributions are weakly convergent and for some constants
$\alpha$, $\beta$ and $K$ such that
$$
\beta\in (0,1]\quad\textrm{and}\quad\alpha\beta>\frac 2{\log_2\left(\frac{4d}{4d-3}\right)}
$$
and
\begin{equation}\label{cond_tightness_Jana}
 \mu\left\{|Y_n(t)-Y_n(s)|\geqslant\varepsilon\right\}\leqslant \frac K{\varepsilon^\alpha
 }\lVert s-t\rVert^{\alpha\beta}
\end{equation}
for any $s$ and $t$ in $[0,1]^d$, any $\varepsilon>0$ and any positive integer $n$ then $(Y_n(\cdot))_{n\geqslant 1}$ converges weakly to some process in
$\mathbb{H}_\gamma([0,1]^d)$ where $0<\gamma <\beta-m/\alpha$. Since the finite-dimensional laws of the process $\{n^{-d/2}S_n(t)\,;\,t\in[0,1^d]\}_{n\geqslant 1}$ are weakly convergent (cf. Theorem \ref{IP_OMD_bis}), it suffices to convert the moment inequality given by Proposition~2 into
an inequality involving $\mu\left\{|S_n(t)-S_n(s)|\geqslant n^{d/2}
\varepsilon\right\}$ in order to check that condition
\eqref{cond_tightness_Jana} is satisfied with $\alpha=p$, $\beta=1/2$ and
$Y_n(t)=n^{-d/2}S_n(t)$. We shall do the proof for $d=2$. Let $s=(s_1,s_2)$ and
$t=(t_1,t_2)$ be fixed in $[0,1]^2$ and $n$ be a positive integer. Without loss of generality, we assume that $s_1>t_1$ and
$s_2<t_2$ (similar arguments can be used to threat the general case). Let $s'_1=k_1/n$ and $t'_1=(l_1+1)/n$ where $(k_1,l_1)$ is the unique element of $\langle n\rangle^2$ such that $k_1/n\leqslant s_1<(k_1+1)/n$ and $l_1/n\leqslant t_1<(l_1+1)/n$. In other words, keeping in mind that $[.]$ denotes the integer part function, we have $s'_1=[ns_1]/n$ and $t'_1=([nt_1]+1)/n$ and similarly, we define $s'_2=([ns_2]+1)/n$ and $t'_2=[nt_2]/n$. With these notations, we have
\begin{align*}
 |S_n(t)-S_n(s)|&=|S_n(t_1,t_2)-S_n(s_1,s_2)|\\
 &\leqslant  |S_n(t_1,t_2)-S_n(t_1,t'_2)|+|S_n(t'_1,t'_2)-S_n(t_1,t'_2)|\\
 &\quad+|S_n(t'_1,t'_2)-S_n(s'_1,s'_2)|+|S_n(s'_1,s'_2)-S_n(s'_1,s_2)|\\
 &\quad+|S_n(s'_1,s_2)-S_n(s_1,s_2)|.
\end{align*}
Since
$$|S_n(t_1,t_2)-S_n(t_1,t'_2)|=(t_2-t'_2)\left|\sum_{i=1}^{[nt_1]}
X_{i,[nt_2]}+(t'_1-t_1)X_{[nt_1]+1,[nt_2]}\right|
$$
and $t_2-t'_2\leqslant 1/n$, we have
\begin{equation}\label{Moment_inequality_for_t_tprime}
\mathbb E\left|S_n(t_1,t_2)-S_n(t_1,t'_2)\right|^p\leqslant \kappa (t_2-t'_2)^pn^{p/2}\mathbb E|X_{0,0}|^p\leqslant \kappa(t_2-t'_2)^{p/2}
\mathbb E|X_{0,0}|^p.
\end{equation}
Similarly,
\begin{align}
\mathbb E\left|S_n(t'_1,t'_2)-S_n(t_1,t'_2)\right|^p
&\leqslant \kappa(t'_1-t_1)^{p/2}\mathbb E|X_{0,0}|^p,\label{Moment_inequality_for_t_tprime_2}\\
\mathbb E\left|S_n(s'_1,s'_2)-S_n(s'_1,s_2)\right|^p
&\leqslant \kappa(s'_2-s_2)^{p/2}\mathbb E|X_{0,0}|^p,\label{Moment_inequality_for_t_tprime_3}\\
\mathbb E\left|S_n(s'_1,s_2)-S_n(s_1,s_2)\right|^p
&\leqslant \kappa(s_1-s'_1)^{p/2}\mathbb E|X_{0,0}|^p\label{Moment_inequality_for_t_tprime_4}.
\end{align}
Moreover, from Proposition \ref{Moment_inequality_under_proj_condition}, for any positive $n$ and any $i$ and $j$ in $\langle n\rangle^2$, we have\begin{equation}\label{Moment_inequality_for_i_over_n_and_j_over_n}
\mathbb{E}\left\vert \frac{1}{n}S_n\left(\frac{i}{n}\right)-\frac{1}{n} S_n\left(\frac{j}{n}\right)\right\vert^p\leqslant \kappa\mathbb{E}\vert X_{0,0}\vert^p\norm{\frac{i}{n}-\frac{j}{n}}^{p/2}.
\end{equation}
Combining \eqref{Moment_inequality_for_t_tprime}, \eqref{Moment_inequality_for_t_tprime_2}, \eqref{Moment_inequality_for_t_tprime_3}, \eqref{Moment_inequality_for_t_tprime_4} and \eqref{Moment_inequality_for_i_over_n_and_j_over_n} and using the elementary convexity inequality
$(a_1+a_2+a_3+a_4+a_5)^p\leqslant 5^{p-1}(a_1^p+a_2^p+a_3^p+a_4^p+a_5^p)$for any non-negative $a_1,a_2,a_3,a_4$ and $a_5$, we derive that
$$
\mathbb E|S_n(t)-S_n(s)|^p\leqslant\kappa\lVert s-t\rVert^{p/2}.
$$
Finaly, using Markov's inequality, we obtain \eqref{cond_tightness_Jana}. The proof of Theorem \ref{IP_Holder} is complete.\\
\\
\textbf{Acknowledgements.} The authors would like to thank D. Voln\'y and Y. Wang for many usefull comments and discussions.

\end{document}